\documentclass{gtart}

\def\ifplaintex{\expandafter\ifx\csname documentclass\endcsname\relax}

\ifplaintex 
\else
\fi

\def\gtm{{\mathsurround=0pt\it $\cal G\mskip-2mu$eometry \&\ 
$\cal T\!\!$opology $\cal M\mskip-1mu$onographs}}    

\def\gtp{{\mathsurround=0pt\it $\cal G\mskip-2mu$eometry \&\ 
$\cal T\!\!$opology $\cal P\!$ublications}}  

\def\recd{{\small Received:\qua\receiveddate\ifx\reviseddate\relax
\else\qquad Revised:\qua\reviseddate\fi\par}} 

\def\volumenumber#1{\def\thevolumenumber{#1}}
\def\volumeyear#1{\def\thevolumeyear{#1}}
\def\volumename#1{\def\thevolumename{#1}}
\def\papernumber#1{\def\thepapernumber{#1}}
\def\pagenumbers#1#2{\def\startpage{#1}\def\finishpage{#2}}
\def\published#1{\def\publishdate{#1}}
\def\received#1{\def\receiveddate{#1}}
\def\revised#1{\def\reviseddate{#1}}
\def\accepted#1{\def\accepteddate{#1}}

\long\def\asciiabstract#1{\long\def\theasciiabstract{#1}}

\let\\\par
\let\thevolumenumber\relax\let\thepapernumber\relax
\let\thevolumeyear\relax\let\startpage\relax
\let\finishpage\relax\let\publishdate\relax\let\receiveddate\relax
\let\reviseddate\relax\let\accepteddate\relax\let\theasciititle\relax
\let\theasciiauthors\relax
\let\theasciiabstract\relax

\let\theerratum\relax\let\theasciiemail\relax
\let\theshortauthors\relax\let\theshorttitle\relax

\def\startpage{1}\def\finishpage{15}\def\thepapernumber{77}

\volumenumber{2}
\volumename{Proceedings of the Kirbyfest}
\volumeyear{1999}

\long\def\maketitlep{   

\count0=\startpage

\gtm\nl        
{\small Volume \thevolumenumber: \thevolumename\nl 
\ifx\theerratum\relax\else Erratum \erratumnumber\nl\fi
Pages \startpage--\finishpage\nl}

\vglue 0.1truein   

{\parskip=0pt\leftskip 0pt plus 1fil\def\\{\par\smallskip}{\ifplaintex\large
\else\Large\fi\bf\thetitle}\par\medskip}   
\vglue 0.05truein 

{\parskip=0pt\leftskip 0pt plus 1fil\def\\{\par}{\sc\theauthors}
\par\medskip}%
 
\vglue 0.03truein 

{\small\leftskip 25pt\rightskip 25pt{\bf Abstract}\stdspace\theabstract

{\bf AMS Classification}\stdspace\theprimaryclass
\ifx\thesecondaryclass\relax\else; \thesecondaryclass\fi\par
{\bf Keywords}\stdspace \thekeywords\par}\vglue 7pt

}   

\font\phead=cmsl9 scaled 950
\font=cmsl9 scaled 1050
\font\pnum=cmbx10 scaled 913
\font\lnum=cmbx10 
\font\pfoot=cmsl9 scaled 950
\font=cmsl9 scaled 1050
\ifplaintex
\headline{\vbox to 0pt{\vskip -4.5mm\line{\small\phead\ifnum
\count0=\startpage ISSN 1464-8997 (on line)
1464-8989 (printed) \hfill {\pnum\folio}\else\ifodd\count0\def\\{ }%
\ifx\theshorttitle\relax\thetitle\else\theshorttitle\fi\hfill{\pnum\folio}
\else\def\\{ and }{\pnum\folio}\hfill\ifx\theshortauthors\relax\theauthors
\else\theshortauthors\fi\fi\fi}\vss}}
\footline{\vbox to 0pt{\vglue 0mm\line{\small\pfoot\ifnum\count0=\startpage
Published \publishdate:\qua\copyright\ \gtp\hfill\else
\gtm, Volume \thevolumenumber\ (\thevolumeyear)\hfill\fi}\vss
}}
\else
\makeatletter
\def\@oddhead{{\small\lhead\ifnum\count0=\startpage ISSN 1464-8997 (on line)
1464-8989 (printed) \hfill {\lnum\number\count0}\else\ifodd\count0
\def\\{ }\ifx\theshorttitle\relax \thetitle \else\theshorttitle\fi\hfill
{\lnum\number\count0}\else\def\\{ and }{\lnum\number\count0}
\hfill\ifx\theshortauthors\relax 
\theauthors\else\theshortauthors\fi\fi\fi}}\def\@evenhead{@oddhead}
\def\@oddfoot{\small\lfoot\ifnum\count0=\startpage Published \publishdate:\qua\copyright\ \gtp\hfill\else
\gtm, Volume \thevolumenumber\ (\thevolumeyear)\hfill\fi}
\def\@evenfoot{@oddfoot}
\makeatother
\fi

\let\maketitlepage\maketitlep
\let\makeshorttitle\maketitlepage
\let\maketitle\maketitlepage

\newwrite\gtoutfile
\long\gdef\makeheadfile{  
{\def\\{, }\def\s{ }
\immediate\openout\gtoutfile head.xxx
\immediate\write\gtoutfile{To: math@arxiv.org}
\immediate\write\gtoutfile{Subject: put OR rep NNNNN:ppppp}
\immediate\write\gtoutfile{--text follows this line--}
\immediate\write\gtoutfile{Proxy-for: \ifx\theasciiauthors\relax
\theauthors\else\theasciiauthors\fi\s<\ifx\theasciiemail\relax\theemail\else\theasciiemail\fi>}
\immediate\write\gtoutfile{\noexpand\\}
\immediate\write\gtoutfile{Authors: \ifx\theasciiauthors\relax
\theauthors\else\theasciiauthors\fi}
{\def\\{ }\immediate\write\gtoutfile{Title: \ifx\theasciititle\relax
\thetitle\else\theasciititle\fi}}
\immediate\write\gtoutfile{Subj-class: GT or SG, GR etc}
\immediate\write\gtoutfile{MSC-class: \theprimaryclass\ifx\thesecondaryclass\relax\else, \thesecondaryclass\fi}
\immediate\write\gtoutfile{Journal-ref: Geom. Topol. Monogr. \thevolumenumber\s
(\thevolumeyear) \startpage-\finishpage}
\immediate\write\gtoutfile{Comments: Published by Geometry and Topology Monographs at}
\immediate\write\gtoutfile{\s\s\s  http://www.maths.warwick.ac.uk/gt/GTMon\thevolumenumber/paper\thepapernumber.abs.html}
\immediate\write\gtoutfile{\noexpand\\}
\immediate\write\gtoutfile{}
\ifx\theasciiabstract\relax
\immediate\write\gtoutfile{\theabstract}\else
\immediate\write\gtoutfile{\theasciiabstract}\fi
\immediate\write\gtoutfile{}
\immediate\write\gtoutfile{\noexpand\\}
\immediate\write\gtoutfile{}
\immediate\closeout\gtoutfile}}  

\def\maketitlepage{\maketitlep\makeheadfile}
\let\makeshorttitle\maketitlepage
\let\maketitle\maketitlepage

\volumenumber{4}
\volumename{Invariants of knots and 3-manifolds (Kyoto 2001)}
\volumeyear{2002}
\papernumber{15}
\pagenumbers{235}{244}
\received{31 October 2001}
\revised{15 May 2002}
\accepted{22 July 2002}
\published{13 October 2002}

\usepackage[rokicki, hideboxes]{boxedeps}

\newtheorem{theorem}{Theorem}
\newtheorem{corollary}[theorem]{Corollary}
\newtheorem{lemma}[theorem]{Lemma}
\newtheorem*{conjecture}{Conjecture}
\theoremstyle{definition}
\newtheorem*{remark}{Remark}

\newtheorem*{Def}{Definition}

\newcommand{\pic}[2]{\BoxedEPSF{#1 scaled #2}}

\def\CA{{\cal A}}
\def\CC{{\cal C}}
\def\CS{{\cal S}}

\def\Xor{\pic{Xor.ART} {400}}
\def\Yor{\pic{Yor.ART} {400}}
\def\Ior{\pic{Ior.ART} {400}}

\def\Idor{\pic{Idor.ART} {400}}
\def\Rcurlor{\pic{Rcurlor.ART} {400}}

\def\Murphy{\pic{Murphy.ART} {400}}
\def\Murphysum{\pic{Murphys.ART} {400}}

\def\Gcross{\pic{Gcross.ART} {300}}
\def\Apos{\pic{Gcross.ART} {250}}
\def\Amix{\pic{Amix.ART} {250}}
\def\Aneg{\pic{Aneg.ART} {250}}
\def\Aij{\pic{Aij.ART} {300}}

\def\Tjrev{\pic{Tjrev.ART} {300}}
\def\Unp{\pic{Unp.ART} {300}}

\def\XC{\pic{XC.ART} {400}}
\def\XCYC{\pic{XCYC.ART} {400}}
\def\YCXC{\pic{YCXC.ART} {400}}
\def\XYC{\pic{XYC.ART} {400}}
\def\TC{\pic{CloseT.ART} {400}}
\def\figthree{\pic{Fnp.ART} {400}}

\def\PsiD{\pic{PsiD.ART} {400}}
\def\PsiX{\pic{PsinpX.ART} {400}}

\def\PsiXY{\pic{PsinpXY.ART} {400}}
\def\PsiXPsiY{\pic{PsinpXPs.ART} {400}}

\newcommand\ds{\displaystyle}
\newcommand\eval[1]{\langle#1\rangle}

\begin{document}

\title{Power sums and Homfly skein theory}                    
\authors{Hugh R. Morton}                  
\address{Department of Mathematical Sciences, University of 
Liverpool\\Peach St, Liverpool, L69 7ZL, UK}                  
\email{morton@liv.ac.uk}
\url{http://www.liv.ac.uk/\char'176su14}

\begin{abstract}   
The Murphy operators in the Hecke algebra
$H_n$ of type $A$ are explicit commuting
elements, whose 
symmetric functions are  central in $H_n$.
In \cite{Murphy} I defined geometrically a  homomorphism
from the Homfly skein 
$\CC$ of the annulus to the centre of each algebra $H_n$, and
found an element $P_m$ in $\CC$, independent of $n$, 
whose image, up to an explicit linear combination with the identity of $H_n$,  is the $m$th power sum of the
Murphy operators. The aim of this paper is to give simple geometric representatives for the elements $P_m$, and to discuss their role in a similar construction for central elements of an extended family of algebras $H_{n,p}$.
\end{abstract}

\asciiabstract{   
The Murphy operators in the Hecke algebra H_n of type A are explicit
commuting elements, whose symmetric functions are central in H_n.  In
[Skein theory and the Murphy operators, J. Knot Theory Ramif. 11
(2002), 475-492] I defined geometrically a homomorphism from the
Homfly skein C of the annulus to the centre of each algebra H_n, and
found an element P_m in C, independent of n, whose image, up to an
explicit linear combination with the identity of H_n, is the m-th
power sum of the Murphy operators. The aim of this paper is to give
simple geometric representatives for the elements P_m, and to discuss
their role in a similar construction for central elements of an
extended family of algebras H_{n,p}.}


\makeshorttitle  

\section*{Introduction}

This article is a combination of material presented at  a conference in
Siegen in January 2001 and at the Kyoto low-dimensional topology workshop in
September 2001. I am grateful to the organisers of both these meetings for
the opportunity to discuss this work.
 I must also thank the Kyoto organisers and RIMS for their hospitality during
the production of this account.

Since much of the material in the talks is already contained in 
\cite{Murphy} I  refer readers there for many of the details, while giving a
brief review here. In this paper I  concentrate on some of the related
results and developments which did not appear in \cite{Murphy}. In particular
I  include a direct skein theory proof in theorem \ref{braidsum} of a result
originally due to Aiston \cite{Adams} giving an unexpectedly simple geometric
representation for a sequence of elements $\{P_m\}$  in the skein of the
annulus  as  sums of a small number of closed braids. The  elements $P_m$
were interpreted in \cite{Adams} as the result of applying Adams operations
to the core curve of the annulus, when the skein of the annulus is viewed as
the representation ring of $sl(N)$ for\eject 
large $N$.  The same elements $P_m$ are shown in \cite{Murphy} to
give rise to the power sums of the Murphy operators in the centre of
the Hecke algebras $H_n$.

I extend the ideas to algebras $H_{n,p}$ based on tangles with strings
allowed in both directions at top and bottom, where a similar uniform
description for much, in fact probably all,  of the centre, can be given.
This suggests analogues of the Murphy operators in these cases, which now
come in two sets $\{T(j)\}, j=1,\ldots,n$ and $\{U(k)\}, k=1,\ldots,p$, of
commuting elements, and provide central elements as the supersymmetric
functions of these sets of elements in the algebra $H_{n,p}$.
 
\section{Framed Homfly skeins}

The skein theory setting  described more fully in \cite{Murphy} is based on
the framed Homfly skein relations,
\begin{eqnarray*}
\Xor\  -\ \Yor& =&(s-s^{-1})\ \Ior\\
\mbox{and }\  \Rcurlor &=&v^{-1}\ \Idor.
\end{eqnarray*}
The  framed Homfly skein ${\cal S}(F)$ of a planar surface $F$, with some
designated input and output boundary points, is defined as linear
combinations of oriented tangles in $F$, modulo Reidemeister moves II and III
and these two local skein relations.
 The coefficient ring can be taken as $\Lambda={\bf Z}[v^{\pm1}, s^{\pm1}]$
with powers of $s^k-s^{-k}$ in the denominators.

When  $F$ is a rectangle with $n$ inputs at the bottom and
$n$ outputs at the top the skein $R_n^n(v,s)={\cal S}(F)$ of
$n$-tangles is an algebra under composition of tangles, isomorphic to the
Hecke algebra
$H_n(z)$, with $z=s-s^{-1}$ and coefficients extended to $\Lambda$.
 This algebra has a presentation with generators
$\{\sigma_i\}, i=1,\ldots,n-1$ corresponding to Artin's elementary braids,
which satisfy the braid relations
 and the quadratic relations
$\sigma_i^2=z\sigma_i+1$.
The braids $T(j), j=1,\ldots,n$, shown in figure 1 make up a set of commuting
elements in $H_n$. These are shown by Ram \cite{Ram} to be the Murphy
operators of Dipper and James \cite{Dipper}, up to linear combination with
the identity of $H_n$.
\begin{center}
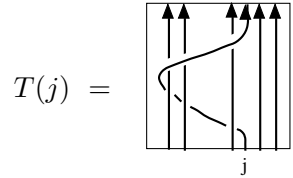
$T(j)\ =\ $ \Murphy\\ {\small Figure 1}
\end{center}
Their properties are discussed further in \cite{Murphy}, where the element
$T^{(n)}$ in figure 2 is shown to represent their sum, again up to linear
combination with the identity.
\begin{center}
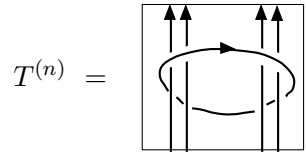
$T^{(n)}\ =\ $ \Murphysum\\ {\small Figure 2}
\end{center}
 In the same paper the power sums $\sum T(j)^m$ are presented as skein
elements which clearly belong to the centre of $H_n$, in terms of an element
$P_m$ in the skein $\CC$ of the annulus. This leads to skein theory
presentations for any symmetric function of the Murphy operators as central
elements of $H_n$, and gives a pictorial view of the result of Dipper and
James which identifies the centre of $H_n$ for generic parameter with these
symmetric functions.

Before giving the  skein presentation for $P_m$ as a sum of $m$ closed braids
in the annulus I discuss briefly the construction of central elements in
$H_n$, and in some extended variants of these algebras.

\begin{Def} Write $H_{n,p}$ for the skein $\CS(F)$ where $F$ is the rectangle
with $n$ outputs and $p$ inputs at the top, and matching inputs and outputs
at the bottom as in figure 3.
\end{Def}
\begin{center}
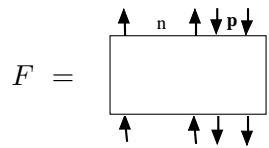
$F\ =\ $ \figthree\\ {\small Figure 3}
\end{center} There is  a natural algebra structure on $H_{n,p}$ induced by
placing oriented tangles one above the other. When $p=0$ we have the Hecke
algebra $H_n=H_{n,0}$. The resulting algebra $H_{n,p}$ has been studied by
Kosuda and Murakami, \cite{Kosuda}, in the context of $sl(N)_q$ endomorphisms
of the module $V^{\otimes n}\otimes \overline{V}^{\otimes p}$, where $V$ is
the fundamental $N$-dimensional module. Hadji \cite{Hadji} has also described
an explicit skein-theoretic basis for it; while there is a linear isomorphism
of $H_{n,p}$ with $H_{(n+p)}$ this is not in general an algebra isomorphism.

\section{The annulus} The skeins $H_{n,p}$ are closely related to the skein
$\CC=\CS(F)$ where $F$ is the annulus. An element  $X\in \CC$, which is
simply a linear combination of diagrams in the annulus, modulo the skein
relations, is indicated schematically as 
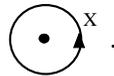
\[\XC\ .\]
The closure map $H_{n,p}\to {\cal C}$, induced by taking a tangle $T$ to its
closure $\hat {T}$ in the annulus, is defined by \[\hat {T}\ =\ \TC\ .\] This
is a $\Lambda$-linear map, whose image we  call ${\CC}_{n,p}$. Every diagram
in the annulus lies in some $\CC_{n,p}$.

The skein $\CC$ has a product induced by placing one annulus outside another, 
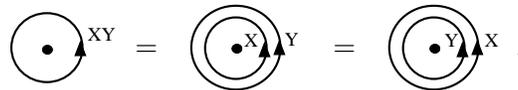
\[\XYC\ =\ \XCYC\ =\ \YCXC .\]
 Under this product $\CC$ becomes a commutative algebra, since diagrams in
the inner annulus can be moved across the outer annulus using Reidemeister
moves II and III.

The evaluation map $\eval{\ }:\CC\to\Lambda$ is defined by setting $\eval X$
to be the framed Homfly polynomial of $X$, regarded as a linear combination
of diagrams in the plane. When the Homfly polynomial is normalised to take
the value 1 on the empty diagram, and hence the value
$\delta=\ds\frac{v^{-1}-v}{s-s^{-1}}$ on the zero-framed unknot, the map
$\eval{\ }$ is a multiplicative homomorphism.

\section{Central elements of $H_{n,p}$} There is  an easily defined algebra 
homomorphism $\psi_{n,p}$ from
${\cal C}
$ to the centre of each algebra $H_{n,p}$, induced from
$D=\PsiD$  by placing $X\in \CC$  around the circle and the identity of
$H_{n,p}$ on the arc, to get
\[\psi_{n,p}(X)\ =\ \PsiX\ \in H_{n,p}.\]
 It is clear that 
\[\psi_{n,p}(XY)\ =\ \PsiXY\ =\ \PsiXPsiY\ =\psi_{n,p}(X)\psi_{n,p}(Y),\] and
that
 the elements $\psi_{n,p}(X)$  all lie in the centre of $H_{n,p}$. It is
shown in \cite{Murphy} that the image of $\psi_{n,0}$, called $\psi_n$ there,
consists of all symmetric polynomials in the Murphy operators in $H_n$ and
so, by \cite{Dipper}, makes up the whole of the centre of $H_n$ in the
generic case. 

It is natural to suspect that the same is true for $H_{n,p}$.
\begin{conjecture} The image of $\psi_{n,p}$ is the whole centre of $H_{n,p}$.
\end{conjecture}
\begin{remark}  I don't have an
 immediate skein theory argument showing that central elements can always be
written as $\psi(X)$, even in the case of  $H_n$.   I suspect that the best
way would be to show that the restriction of $\psi_{n,p}$ to $\CC_{n,p}$
followed by closure is an isomorphism, and couple this with an upper bound on
the dimension of the centre, although I don't know of an algebraic
determination of the centre in the general case $n,p>0$. 
\end{remark}

In \cite{Murphy} the element $P_m\in \CC$ is constructed, by manipulating
generating functions, as a polynomial in elements $\{h_i\}$, while each $h_i$
is itself initially a linear combination of $i!$ closed braids. It satisfies
\[\psi_{n,0}(P_m-\eval{P_m})=(s^m-s^{-m})v^{-m}\sum_{j=1}^n T(j)^m,\] for
every $n$. This result can be generalised to $H_{n,p}$ as follows.
\begin{theorem} The central element $\psi_{n,p}(P_m)$ of $H_{n,p}$ can be
written, up to linear combination with the identity, as the power sum
difference \[v^{-m}\sum_{j=1}^n T(j)^m-v^m\sum_{k=1}^pU(k)^m,\] for some
commuting elements $\{T(j)\}\cup\{U(k)\}$.
\end{theorem}

\begin{proof} Apply the arguments of \cite{Murphy} using the skein $\CA$ of
the annulus with a boundary input and output to show that
\[\psi_{n,p}(P_m)-\eval{P_m}{\rm Id}=(s^m-s^{-m})\left(v^{-m}\sum_{j=1}^n
T(j)^m-v^m\sum_{k=1}^p U(k)^m\right),\]  where 
\[T(j)\ =\ \Tjrev\ ,\quad U(k)\ =\ \Unp\ .\]
\end{proof} 
\begin{remark}
Every central element of the form $\psi_{n,p}(X)$ is a
polynomial in the elements $\psi_{n,p}(P_m)$. Stembridge 
\cite{Stembridge} defines a supersymmetric polynomial in two sets of
commuting variables $\{x_i\}$ and $\{y_j\}$ to be a polynomial which is
symmetric in $\{x_i\}$ and $\{y_j\}$ separately, and which is independent of
$t$ after making the substitution $x_1=y_1=t$. He then characterises them as
polynomials in the power sum differences $\sum x_i^m-\sum y_j^m$.  The
central elements
$\psi_{n,p}(X)$ can thus be described  as  supersymmetric polynomials in the
sets of commuting elements $\{v^{-1}T(j)\}$ and $\{vU(k)\}$.

These elements  might be
considered
  as a set of Murphy operators in the  general setting.  There is certainly
some measure of choice for the set of elements; conjugating all of them by a
fixed element will not alter their supersymmetric polynomials, while further
choices are also possible. 
\end{remark}

\section{Geometric representatives of $P_m$}

The main result of this paper is to give a geometrically simple 
representative for $P_m$ as a sum of  $m$ closed braids.

Write $A_m \in {\cal C}$ for the closure of the
$m$-braid \[\sigma_{m-1}\cdots \sigma_2\sigma_1\ =\ \Gcross\ .\] The central
element $T^{(n)}\in H_n$ can be written as $T^{(n)}=\psi_{n,0}(A_1)$, so we
can take
 $P_1$ to be the core curve $A_1$ of the annulus.

The construction of the elements $P_m$ in \cite{Murphy} makes use of the
sequence of elements $h_i\in \CC_{i,0}$ which arise as the closure of one of
the two most basic idempotents in $H_i$. The properties of these elements
which will be used here are their relations to the elements $A_m$. These are
expressed compactly in terms of formal power series with coefficients in the
ring $\CC$ in theorem \ref{AH}. 

\begin{theorem}[Morton] \label{AH}
 \[A(t)=\frac{H(st)}{H(s^{-1}t)},\] where
\[H(t)=1+\sum_{n=1}^\infty h_n t^n\] and
 \[A(t)=1+(s-s^{-1})\sum_{m=1}^\infty A_mt^m.\]
\end{theorem} 
\begin{proof} This is given in \cite{Murphy}, using simple skein properties
of the elements $h_n$ and the skein $\CA$ of the annulus with a single input
on one boundary curve and a matching output on the other.
\end{proof}

The elements $P_m$ are defined in \cite{Murphy} by the formula
\[\sum_{m=1}^\infty \frac{P_m}{m}t^m=\ln(H(t)).\]
Now every skein admits a {\em mirror map},
$\overline{\phantom{w}}:\CS(F)\to\CS(F)$ induced by switching all crossings
in a tangle, coupled with inverting $v$ and $s$ in $\Lambda$. When this is
applied to the series $A(t)$ the result is the series
\[\overline{A}(t)=1-z\sum_{m=1}^\infty \overline{A}_mt^m,\] where
$\overline{A}_m$ is the mirror image of the closed braid $A_m$, and
$z=s-s^{-1}$.

 The series $H(t)$ is invariant under the mirror map, as shown in
\cite{Murphy}, and hence so is $P_m$, although I will not need  these facts
here.  The simple geometric representative for $P_m$ in theorem
\ref{braidsum} can be thought of as `almost' $A_m$, but averaged in  a way
that ensures mirror symmetry.
\begin{Def} Write $A_{i,j}$ for the closure of the braid 
\[\sigma_{i+j}^{-1}\sigma_{i+j-1}^{-1}\ldots\sigma_{i+1}^{-1}\sigma_i\ldots\sigma_2\sigma_1
=\ \Aij\] on $i+j+1$ strings.
\end{Def} Then  $A_{i,j}$ is  the closed braid with $i$ positive and $j$
negative crossings which results from switching the last $j$ crossings of
$A_{i+j+1}$, with $A_m=A_{m-1,0}$ and $\overline{A}_m=A_{0,m-1}$.

\begin{theorem}[Aiston] \label{braidsum} The element $P_m$ is a multiple of
the sum of all $m$ of the closed braids $A_{i,j}$ which have $m$ strings.
Explicitly,
\[[m]P_m=\sum_{{i=0}\atop{i+j=m}}^{m-1}A_{i,j-1}=\overline{A}_m+\ldots+{A}_m
=\Aneg +\cdots +\Amix +\cdots +\Apos,\]  where $\ds
[m]=\frac{s^m-s^{-m}}{s-s^{-1}}$.
\end{theorem}

While Aiston's proof required a diversion through results about $sl(N)_q$
representations, the proof here uses only the skein-based result of theorem
\ref{AH} and one further simple skein theory lemma.

\begin{lemma}\label{Adiff} For all $i,j>0$ we have
\[A_{i,j-1}-A_{i-1,j}=zA_{i-1,0}A_{0,j-1}\ (=zA_i\overline{A}_j).\]
\end{lemma}
\begin{proof} Apply the skein relation to the $i$th positive crossing in
$A_{i,j-1}$. Switching this crossing gives $A_{i-1,j}$ while smoothing it 
gives the product $A_i\overline{A}_j$.
\end{proof}

\begin{corollary}\label{mirror}
\[A(t)\overline{A}(t)=1.\]
\end{corollary} 
\begin{proof} We have
\[
z\sum_{{i=1}\atop{i+j=m}}^{m-1}A_i\overline{A}_j=A_{m-1,0}-A_{0,m-1}=A_m-\overline{A}_m,m>1,\]
by lemma \ref{Adiff}. The coefficient of $t^m$ in $A(t)\overline{A}(t)$ is
\[zA_m-z\overline{A}_m-z^2\sum_{i+j=m}A_i\overline{A}_j =0, m>0.\]
\end{proof}
\begin{remark} This  result appears in \cite{Murphy}; the proof here  does
not need theorem  \ref{AH} or properties of the mirror map.
\end{remark}

\begin{proof} [Proof of theorem \ref{braidsum}] Write
\[\Pi_m=\sum_{{i=0}\atop{i+j=m}}^{m-1}A_{i,j-1}.\] Then
\[\sum_{{i=1}\atop{i+j=m}}^{m-1}i(A_{i,j-1}-A_{i-1,j})+\Pi_m=mA_m,\] giving
\[\Pi_m=mA_m-z\sum_{{i=1}\atop{i+j=m}}^{m-1}iA_i\overline{A}_j.\] This is the
coefficient of $t^m$ in
$\left(\sum^\infty_{i=1}iA_it^i\right)\overline{A}(t)$. Now \[z
\sum^\infty_{i=1}iA_it^i=t\frac{d}{dt}(A(t)).\] We then have
\begin{eqnarray}
z\sum_{m=1}^\infty\Pi_mt^{m-1}&=&\frac{d}{dt}(A(t))\overline{A}(t)\nonumber\\
&=&\frac{d}{dt}\ln(A(t)),\label{log}
\end{eqnarray} since $\overline{A}(t)=A(t)^{-1}$ by corollary \ref{mirror}.
Now \begin{eqnarray} \ln A(t)&=& \ln H(st)-\ln H(s^{-1}t),\nonumber\\ &=&
\sum_{m=1}^\infty (s^m-s^{-m})\frac{P_m}{m}t^m,\label{Pm}
\end{eqnarray} by theorem \ref{AH} and the definition of $P_m$. Comparing the
coefficients of $t^{m-1}$ in (\ref{log}), using(\ref{Pm}),  gives
\[(s-s^{-1})\Pi_m=(s^m-s^{-m})P_m,\] and hence the result.
\end{proof}

\Addresses\recd
\end{document}